\documentclass[12pt,amstex]{article}

\usepackage{amsmath}
\usepackage{amsfonts}

\usepackage{amsthm}

\newtheorem{theo}{Theorem}
\newtheorem{coro}{Corollary}

\newtheorem{prop}{Proposition}
\newtheorem{lemm}{Lemma}

\newcommand{\lbl}{\label}

\newcommand{\eq}[1]{$(\ref{#1})$}

\def\Bi{{\rm Bin}}
\def\rng{{\rm rng}}
\def\Mult{{\rm Mult}}
\def\Bin{{\rm Bin}}
\def\N{\mathbb{N}}
\def\E{\mathbb{E}}

\def\0{{\bf 0}}
\def\bx{{\bf x}}
\def\bX{{\bf X}}

\def\bp{{\bf p}}
\def\phat{\hat{p}}
\def\bphat{\hat{{\bf p}}}
\def\N{\mathbb{N}}

\def\R{\mathbb{R}}

\def\G{{\cal G}}

\renewcommand{\E}{\mathbb E \,}

\newcommand{\eqco}{\setcounter{equation}{0}}
\newcommand{\thco}{\setcounter{theo}{0}}
\newcommand{\prco}{\setcounter{prop}{0}}
\newcommand{\laco}{\setcounter{lemm}{0}}
\newcommand{\coco}{\setcounter{coro}{0}}
\newcommand{\cjco}{\setcounter{conj}{0}}

\newcommand{\deco}{\setcounter{defn}{0}}
\newcommand{\allco}{\eqco  \thco \prco \laco \coco \cjco \deco}
\newcommand{\LL}{{\cal L}}

\newcommand{\Cov}{{\rm Cov}}

\newcommand{\Var}{{\rm Var}}

\newcommand{\var}{{\rm Var}}

\newcommand{\XX}{{\bf X}}

\newcommand{\BB}{{\cal B}}

\newcommand{\F}{{\cal F}}
\renewcommand{\G}{{\cal G}}

\newcommand{\tih}{{\tilde{h}}}

\newcommand{\eps}{\varepsilon}
\def\bdm{\begin{displaymath}}
\newcommand{\edm}{\end{displaymath}}
\def\benu{\begin{enumerate}}
\def\eenu{\end{enumerate}}
\def\beqn{\begin{equation}}
\def\eeqn{\end{equation}}
\def\be{\begin{equation}}
\def\ee{\end{equation}}
\def\bea{\begin{eqnarray}}
\def\eea{\end{eqnarray}}
\newcommand{\bean}{\begin{eqnarray*}}
\newcommand{\eean}{\end{eqnarray*}}
\newcommand{\bear}{\begin{eqnarray}}
\newcommand{\eear}{\end{eqnarray}}

\renewcommand{\epsilon}{\varepsilon}

\newcommand{\BE}{Berry-Ess\'een~}

\def\R{\mathbb{R}}

\def\qed{\hfill\hbox{${\vcenter{\vbox{
    \hrule height 0.4pt\hbox{\vrule width 0.4pt height 6pt
    \kern5pt\vrule width 0.4pt}\hrule height 0.4pt}}}$}}

\def\Sbar{{\bar S}}

\begin{document}
\title{\bf Normal approximation for isolated balls in an urn
allocation model
}

\author{Mathew D. Penrose$^{1,2}$\\
\\
{\normalsize{\em
 University of Bath}} }
\maketitle

 \footnotetext{ $~^1$ Department of
Mathematical Sciences, University of Bath, Bath BA1 7AY, United
Kingdom: {\texttt m.d.penrose@bath.ac.uk} }
 \footnotetext{ $~^2$
 Partly supported by the Alexander
 von Humboldt Foundation
though a Friedrich Wilhelm  Bessel Research Award.
}

\maketitle

%

\begin{abstract}
Consider throwing $n$ balls at random into $m$ urns, each
ball landing in urn $i$ with probability $p_i$.  Let $S$ be the
 resulting number of singletons, i.e., urns containing just one ball.
We give an error bound for the Kolmogorov distance from $S$ to the normal,
 and estimates on its variance.  These  show that if  $n$, $m$ and
 $(p_i, 1 \leq i \leq m)$ vary in such a way that $\sup_i p_i = 
O(n^{-1})$, then $S$ satisfies a CLT if and only if $n^2 \sum_i p_i^2$
tends to infinity,  and demonstrate an optimal rate of convergence in
 the CLT in this case.  In the uniform case $(p_i \equiv m^{-1})$ 
with $m$ and $n$ growing proportionately, we provide bounds with better
  asymptotic constants.  The proof of the error bounds are based on
 Stein's method via size-biased couplings.
\end{abstract}

\section{Introduction}
\lbl{secintro}
\allco
Consider the classical occupancy scheme, in which
each of $n$ balls is placed  independently at random
in one of $m$ urns, with probability $p_i$
of going into the $i$th urn ($p_1 + p_2 + \cdots + p_m=1$).
 If $N_i$ denotes the number of
balls placed in  the $i$th urn, then $(N_1,\ldots,N_m)$
has the multinomial distribution  $\Mult(n; p_1,p_2,\ldots, p_m).  $
A special case of interest is the so-called uniform case
where all the $p_i$ are equal to $1/m$.

A  much-studied quantity is the number of occupied
urns, i.e. the sum $\sum_i {\bf 1} \{N_i > 0 \}$.
This quantity, scaled and centred, is  known to be
asymptotically normal as $n \to \infty$ in the uniform
case with 
$m \propto n$,
and a
Berry-Ess\'een  bound for the discrepancy from the normal,
tending to zero at the optimum rate, was obtained
for the uniform case 
by Englund \cite{Eng}, 
and for the general (nonuniform) case,
 with a less explicit error bound,
 by Quine and Robinson \cite{QR}.
More recently, Hwang and Janson \cite{HJ} have obtained a local limit theorem.
A variety of applications are mentioned in \cite{HJ}
(`coupon collector's problem,  species trapping, birthday paradox,
polynomial factorization,
statistical linguistics, memory allocation, statistical
physics, hashing schemes and so on').  Also noteworthy 
are the monographs
 by Johnson and Kotz \cite{JK} and
by Kolchin {\em et al.} \cite{KSS}; the latter is 
mainly concerned
with models of this type, giving results
for a variety of limiting regimes  for the growth 
of $m$ with $n$ (in the uniform case) and also 
in some of the non-uniform
cases. There has also been recent interest in the
case of infinitely many urns with the probabilities
$p_i$ independent of $n$ \cite{BarGne,GHP}. 

In this paper we consider the number of isolated
balls, that is, the sum $\sum_i {\bf 1} \{N_i = 1\}$.
This quantity seems just as natural  
an object of study as the number of occupied  urns,
if one thinks of the model in terms of the balls
rather than in terms of the urns. For example, in
the well-known birthday paradox, this quantity represents
the number of individuals in the group who have a unique
birthday.  

In the  uniform  case, we obtain  
 an explicit Berry-Ess\'een  bound for the discrepancy of
the number  of isolated balls  from the normal,
tending to zero at the optimum rate when $m \propto n$.
In the non-uniform case we obtain a similar result
with a larger constant, also finding  upper and lower
bounds which show that the variance of the number of isolated balls
 is $\Theta(n^2 \sum_i p_i^2)$. The proof of these bounds, in Section 
\ref{secvarpf}, is based on martingale difference techniques and
somewhat separate from the other arguments in the paper.

Our Berry-Ess\'een results for the number of isolated balls
are analogous to the main results of  
 \cite{Eng} (in the uniform case) and
 \cite{QR} (in the non-uniform case)
for the number of occupied urns.
Our proofs, however, are entirely different.
We adapt a method used recently by Goldstein and Penrose  \cite{GP}
 for a
problem in stochastic 
geometry (Theorem 2.1 of \cite{GP}).


Our method does not involve either characteristic functions,
or first Poissonizing the total number of balls;
in this, it differs from most of the approaches
to problems of this type adopted in the past. As remarked in \cite{HJ} 
`almost all previous approaches rely, explicitly or implicitly,
on the widely used Poissonization technique', and this
remark also applies to \cite{HJ} itself.
One exception is Chatterjee \cite{Chat}, who uses 
a method not involving Poissonization to give an error bound
with the optimal rate of decay (with unspecified constant)
for the Kantorovich-Wasserstein distance
(rather than the Kolmogorov distance, as here)
between the distribution of the number of occupied urns and the normal,
in the uniform case.

We believe that
our approach can be adapted to 
 the number of urns containing $k$ balls,
for arbitrary fixed $k$, but
these might require significant amounts of extra work,
so we restrict ourselves here to the case
with $k=1$.

Our approach
is based on 
 size-biased couplings.
 Given a nonnegative random variable $W$ with
finite mean $\mu=EW$, we say $W'$ has the $W$ size biased
distribution if $P[W' \in dw] = (w/\mu) P[W \in dw]$, or more formally,
if
\bea
\quad
\quad
\E[Wf(W)]=\mu \E f(W') \quad
 \mbox{for bounded continuous functions $f$.}
\lbl{formalsb}
\eea
Lemma \ref{LGthm}  below, due to Goldstein  \cite{Gold},
tells us that if one can find coupled realizations of $W$ and $W'$
which are in some sense close,
then one may be able to
find a good \BE bound for $W$.
It turns out that this can be done for the number of non-isolated balls.

\section{Results}
\label{secgr}
\allco
Let $n \in \N$ and $m = m(n) \in \N$ with $m \geq 4$.
Let $\bp^{(n)} = ( p^{(n)}_x, 1 \leq x \leq m)$
be a probability mass function on $[m] := \{1,2,\ldots,m\}$,
with $p_x^{(n)} > 0 $ for all $x \in [m]$.
Let $X$ and $X_i, 1 \leq i \leq
n$ be independent and identically distributed random variables
with probability mass function  $\bp = \bp^{(n)}$ (we shall
often suppress the superscript $(n)$). Define $Y = Y(n)$  
by 
\bea
M_i := -1 + \sum_{j=1}^n   {\bf 1}
 \{ X_j = X_i\}; ~~~~
Y : = \sum_{i=1}^n     {\bf 1}
 \{  M_i > 0
\}.
\lbl{Yndef}
\eea
In terms of the   urn scheme described in Section \ref{secintro},
the probability
of  landing in Urn $x$ is $p_x$
 for each ball,
$X_i$ represents the location of the $i$th ball,
$M_i$ represents the number of other balls located in
the same urn as the $i$th ball, and
$Y$ represents the number of non-isolated balls, where a ball is
said to be isolated if no other ball is placed in the same urn 
as it is. Thus $n-Y$ is the number of isolated balls,
or in other words, the number of urns which contain a single ball.

Let $Z$ denote a standard normal random variable,
and let $\Phi(t) := P[Z \leq t] = (2 \pi)^{-1/2} \int_{-\infty}^t
\exp(-x^2/2) dx$. Given any random variable $W$ with
finite mean $\mu_W$ and standard deviation  $\sigma_W$
satisfying $0 < \sigma_W < \infty$, define
$$
D_W : = \sup_{t \in \R} \left| P\left[\frac{W- \mu_W}{\sigma_W} \leq t \right]
 - \Phi(t) 
\right|,
$$
 the so-called Kolmogorov
distance between the distribution of $W$ and the normal. 
We are concerned with estimating $D_Y$.

We refer to the  case where $p_x = m^{-1}$
for each $x \in [m]$ as {\em the uniform case}.
 Our main result for the uniform case provides a
 normal approximation error bound for  $Y$, which is
 explicit modulo computation of
 $\mu_Y $ and $\sigma_Y $, and goes as follows.

\begin{theo}
\lbl{thm1}
In the uniform case,
if $\sigma_Y^{3} \geq 24 \mu_Y$, then
\bea
D_Y
\leq
\frac{0.8}{\sigma_Y} + 
 \left( \frac{\mu_Y }{\sigma_Y^{2}} \right)
 \left( \frac{256}{ \sigma_Y} + \frac{32}{ \sigma_Y^{2}}
+ 23
\sqrt{\eta(n,m)}
\right)
\label{thm1eq}
\eea
with
\bea
\eta(n,m) :=
 \frac{16}{n} + \frac{4}{n(n-1)} + \frac{24}{m} 
\left( 2 + \frac{n}{m-3} + \frac{n}{m} \right) 
\lbl{etaYdef}
\eea
\end{theo}
For 
 asymptotics in the uniform case, we
 allow $m = m(n)$ to vary with $n$. We concentrate on the
case  where 
 $m = \Theta(n)$. In this case  
both $\mu_Y$ and $\sigma_Y^2$ turn out to be  $\Theta(n)$ as $n \to \infty$,
 and
thus Theorem \ref{thm1} implies
$D_Y$ 
is $O(n^{-1/2})$  in this regime.
More formally, we have the following.
\begin{theo}
\lbl{coro1}
Suppose $n,m$ both go to infinity in a linked manner, in
such a way that $n/m \to \alpha \in (0,\infty)$.
Then with $g(\alpha) :=
(  e^{-\alpha} 
-  e^{-2 \alpha}( \alpha^2 -  \alpha + 1) )^{1/2}$, 
we have in the uniform case that
$g(\alpha) > 0$ and
\bea
\limsup_{n \to \infty} n^{1/2}   
D_Y
\leq
\frac{0.8}{g(\alpha)} +  256 \left(
\frac{1- e^{-\alpha} }{g(\alpha)^3}
\right)
	 \nonumber \\
+ 92 \left( \frac{1 - e^{-\alpha} }{g(\alpha)^2}  \right)
 (1 + 3 \alpha ( 1 +  \alpha  ) )^{1/2}. 
\lbl{0901d}
\eea
\end{theo}
 In the case $\alpha = 1$, the right hand side of \eq{0901d}, 
rounded up to the nearest integer, comes to 2236.
 Theorems \ref{thm1} and \ref{coro1}  are proved in
 Section \ref{proof-1}.

We now state our results for the general (non-uniform) case. 
Given $n$ we define the parameters
\bea
 \|\bp\| := \sup_{x \in [m]} (  p_x) ; ~~~~~~~~~ 
\gamma = \gamma(n)  := \max( n \|\bp\|, 1).
\lbl{gammadef}
\eea
For the large-$n$ asymptotics we 
essentially assume that $\gamma(n)$ remains bounded,
or at least grows only slowly with $n$; see Corollary \ref{asympthm} below.
First we give a non-asymptotic result.
\begin{theo}
\lbl{nonunithm}
It is the case that
\bea
D_Y  
\geq \min 
\left(  
1/6
, 
( 8 \pi e)^{-1/2} \sigma^{-1} 
\right),
\lbl{thm2eq2}
\eea
and if
\bea
\|\bp\| \leq 1/11
\lbl{sixth}
\eea
and also
\bea
n \geq 83 \gamma^2 (1 + 3\gamma + 3 \gamma^2) e^{1.05 \gamma},
\lbl{nbig1}
\eea 
 then 
\bea
D_Y 
\leq
 8165 \gamma^2 e^{2.1 \gamma}
 ( 577  + 23 C(\gamma)) 
 \sigma_Y^{-1},
\lbl{thm2eq}
\eea
with
\bea
\quad
\quad
\quad
C(\gamma): =
10 (82 \gamma^7 +  82 \gamma^6 +  80 \gamma^5 + 47 \gamma^4 +
 12 \gamma^3 + 12 \gamma^2)^{1/2}.  
\lbl{Cdef}
\eea
\end{theo}
It is of use in proving Theorem \ref{nonunithm}, and also
of independent interest, to estimate the variance $\sigma^2_Y$
in terms of the original parameters $(p_x,x \in [m])$,
and our next result does this.
Throughout, we write $\sum_x$ for $\sum_{x=1}^m$.
\begin{theo}
\lbl{varthm}
It is the case that
\bea
\Var Y \leq 8 n^2 \sum_x  p_x^2,
\lbl{genvarub}
\eea
and if \eq{sixth} and \eq{nbig1}  hold, then 
\bea
\Var Y \geq (7776)^{-1} \gamma^{-2} e^{-2.1 \gamma}  n^2 \sum_x  p_x^2.
\lbl{genvarlb}
\eea
\end{theo}
If 
 $\gamma(n)$ remains bounded, i.e.
$\sup_n \gamma(n) < \infty$,
 then both \eq{sixth} and \eq{nbig1} hold
for large enough $n$.
Hence,
the following asymptotic result
is immediate from Theorems \ref{nonunithm} and \ref{varthm}.
\begin{coro}
\label{asympthm}
Suppose 
$\sup_n \gamma(n) < \infty$.
Then the following three conditions are equivalent:
\begin{itemize}
\item $n^2 \sum_x p_x^2 \to \infty$ as $n \to \infty$ ;
\item 
$\sigma_Y \to \infty $ as $n \to \infty$ ; 
\item
$(Y - \E Y)/\sigma_Y $ converges in distribution to $Z$
as $n \to \infty$.
\end{itemize}
If these conditions hold, then 
$$
D_Y 
 = \Theta(\sigma_Y^{-1}) = \Theta \left( \left(n^2 \sum_x p_x^2 
\right)^{-1/2}
 \right).
$$
\end{coro}

{\em Remarks.}  
In the uniform case, Theorem \ref{coro1} provides an alternative proof 
of the  central limit theorem  for 
 $Y$ when $m = \Theta(n)$
(see
Theorem II.2.4 on page 59 of \cite{KSS}),
 with error bounds 
converging to zero at the optimum rate.  
Corollary \ref{asympthm} shows that in the uniform case, if
$n^2 /m \to \infty$ and $n/m$ remains bounded, then 
$D_Y = \Theta((n^2/m)^{-1/2})$.
Corollary \ref{asympthm} 
overlaps
Theorem III.5.2 on
page 147 of  \cite{KSS}
but is under  weaker conditions than those in  \cite{KSS},
and provides error bounds not given in \cite{KSS}.


The condition that $\gamma(n)$ remain bounded, in 
Corollary \ref{asympthm},
 is also required by 
\cite{QR} for the analogous
\BE type result for the number of occupied boxes,
though not by \cite{HJ}  for the local limit theorem for that quantity. 
In \eq{sixth}, which is used for the non-asymptotic bounds,
the bound of $\frac{1}{11}$ could be replaced by any constant less than 
$\frac{1}{3}$ without changing anything except the constants
in \eq{thm2eq} and \eq{genvarlb}.


As always
(see remarks in \cite{HJ}, \cite{QR2}, \cite{BarGne}),
 it might be possible to obtain similar results to 
those presented here by other methods.
However, to do so
appears to be a non-trivial task.
In \cite{KSS} the count 
of the number of isolated balls is treated
separately, and  differently, from the count of occupied urns or
the count of urns with $k$ balls, $k=0$ or $k \geq 2$.
 Poisson approximation methods might be of use in some
 limiting regimes (see \cite{BHJ}, Chapter 6), but not when
 the ratio between $\E[Y]$ and $\Var[Y]$ remains bounded but is
not asymptotically 1, which is typically the case here.

Exact formulae can be written down for the 
probability mass function and cumulative distribution 
function of $Y$.
For example, using (5.1) on page 99 of \cite{Feller1},
 the cumulative distribution of $Y$ may be written as
\bean
P[ Y \leq n -k ] = P [n-Y \geq k ] = \sum_{j=k}^m (-1)^{j-k}
{j-1 \choose k-1} S_j 
\eean
with  $S_j$ a sum of probabilities that $j$ of the urns contain one ball
each, i.e.
$$
S_j = \sum_{x_1 < x_2 < \cdots < x_j \leq m}
\frac{n!}{(n-j)!}
p_{x_1} \cdots p_{x_j} \left(
 1 - \sum_{i=1}^j p_{x_i} \right)^{n-j} .
$$
We shall not use this formula
in obtaining our normal approximation
results.

\section{Lemmas}
\allco

A key tool in our  proofs is 
 the following result, which is a special case
of Theorem 1.2 of
\cite{Gold}, and is proved there via Stein's method.
\begin{lemm}
\lbl{LGthm} \cite{Gold}
Let $W \ge 0$ be a random variable with mean $\mu$ and variance
$\sigma^2 \in (0,\infty)$, and let $W^s$ be defined on the same
space, with the $W$-size biased distribution. If $|W^s-W| \le B$
for some $B \le \sigma^{3/2}/\sqrt{6 \mu}$, then
\bea
D_W 
 \le
\frac{0.4 B}{\sigma} + \frac{\mu}{\sigma^2}\left(\frac{64 B^2}{\sigma}
+ \frac{4 B^3}{\sigma^2} 
+ 23 \Delta 
\right) 
 \lbl{delbd}
\eea
 where
 \bea \label{def-Delta} \Delta : =
\sqrt{\mbox{\Var}(\E (W^s-W|W))}. \eea
\end{lemm}

Our next lemma is concerned with the construction
of  variables with  size-biased distributions.
\begin{lemm}
\lbl{sblem}
Suppose $W$ is a sum of exchangeable
indicator variables $\xi_1,\ldots,\xi_n$,
with $P[W>0]>0$.
Suppose
$\xi'_1,\ldots,\xi'_n$
are variables
 with joint distribution
$$
{\cal L}(\xi'_1,\ldots,\xi'_n)={\cal L}(\xi_1,\ldots,\xi_n|\xi_1=1).
$$
Then the variable
$
W'=\sum_{i=1}^n \xi'_i
$
has the $W$ size biased distribution.
\end{lemm}
{\em Proof.} See Lemma 3.3 of \cite{GP}.
$\qed$ \\

Let $\Bi(n,p)$ denote the binomial distribution with
parameters $n \in \N$ and $p \in (0,1)$.
The following lemma will be used for constructing the desired close coupling
of our variable of interest $Y$, and its size biased version $Y'$,
so as to be able to use Lemma \ref{LGthm}.
\begin{lemm}
\lbl{bincouplem} Let $\nu \in \N$  and $p \in (0,1)$. Suppose $N
\sim \Bi(\nu,p)$.
Define $\pi_k$ for $k=0,1,2,\ldots,\nu$,
by
\bea
\pi_k :=
\left\{ \begin{array}{lcl}
\frac{
P[N > k |N >0]
 - P[N > k] }{P[N=k](1 - (k/\nu))}
 & \mbox{if} & 0 \leq k \leq \nu -1
\\
0 & \mbox{if} & k= \nu.
\end{array}
\right.
\lbl{pikdef}
\eea
Then $0 \leq \pi_k \leq 1$ for each $k \in \{0,1,\ldots,\nu\}$.
\end{lemm}
{\em Proof.}  See Lemma 3.5 of \cite{GP}.
$\qed$ \\

Our next lemma is a bound on correlations between variables
associated with different balls in the
urn model.
Recall the definition of $M_i$ at \eq{Yndef}

\begin{lemm}
\lbl{momlem}
Let $k \in \N$ with $2 \leq k \leq n$, and
 suppose that for $i=1,\ldots,k$, $\psi_{i}$
 is a real-valued function
defined on $\{0\} \cup [n-1]$, with
 $\E[\psi_{1}(M_1) ]=0$,  and set
$
\|  \psi_{i} \| : = \sup_{\ell \in [n] } \{ |\psi_{i}(\ell -1) |\}
$
and $\rng(\psi_i)
 : = \sup_{\ell \in [n] } \{ \psi_{i}(\ell -1) \} - 
  \inf_{k \in [n] } \{ \psi_{i}(\ell -1) \} $. 
Then in the uniform case,
\bean
\quad 
\left| \E \left[ \prod_{i=1}^k \psi_{i}(M_i) \right] \right|
\leq \frac{k-1}{m}
\left( \prod_{i=2}^k \|\psi_{i} \| \right)
  \rng(\psi_1) 
\left( 
2
+
\left( 
\frac{n}{m-k +1}
\right) 
+
 \frac{n}{m}
 \right).
\eean
\end{lemm}
{\em Proof.}
Set $W:= \prod_{i=1}^k \psi_{i}(M_i) $.
Write $\bX^{(k)}$ for $ (X_1,\ldots,X_k)$ 
 and $\bx$ for $k$-vectors $\bx = (x_1,\ldots,x_k)$  of possible values
of  $\bX^{(k)}$.
Let $F$ be the set of vectors
$\bx = (x_1,\ldots,x_k)$ such that $x_1  
\neq x_j$ for $j=2,\ldots,k$, so that $\{\bX^{(k)} \in F\}$ 
is the
 event that each of Balls $2,\ldots,k$ lands in
a different urn from Ball 1. Then $P[\bX^{(k)} \in F^c] \leq (k-1)/m$, and
$|W|\leq \prod_{i=1}^k \|\psi_i\|$, so that
\bea
|\E W | 
=
\left|
 E [ W|\bX^{(k)} \in F^c] P[\bX^{(k)} \in F^c] 
+ \sum_{\bx \in F} P[\bX^{(k)}= \bx]
 \E [ W|\bX^{(k)} = \bx] 
  \right| 
\nonumber \\
\leq 
((k-1)/m)
\prod_{i=1}^k \|\psi_i\|
+  \sup_{\bx \in F}
 | \E [W |\bX^{(k)} = \bx]|.
\lbl{0422c}
\eea
Fix $\bx \in F$. We group the urns into three `boxes'. 
Let  Box 1 consist of the urn containing Ball 1,
and let Box 2 be the union of the urns containing Balls $2,3,\ldots,k$;
 this could be the union of any number up to $k-1$ of urns depending
on how many of $x_2, x_3, \ldots, x_k$ are distinct,
but  since we assume $\bx \in F$ does not overlap  Box 1.
Let Box 3 consist of all other urns except those
in Box 1 or Box 2.
For $i=1,2,3$,
   let $N_i$ be the number of balls in  Box $i$, other
than Balls $1,\ldots,k$.  Let $h(k)$ be the
expected value of 
$\prod_{i=2}^k\psi_i(M_i) $,
given $\bX^{(k)} = \bx$ and  given that $N_2=k$. Then
\bea
 \E \left[ W | \bX^{(k)} = \bx \right] =
 \E [ \psi_1(N_1) h(N_2) ].
\lbl{0422d} \eea
Also, given $\XX^{(k)} = \bx$,
 $(N_i)_{i=1}^3$ have the multinomial distribution
\bea
\left( N_1,N_2,N_3 \right) \sim \Mult
\left( n-k;
\frac{1}{m},
\frac{a}{m},
1- \frac{1+a}{m} \right),
\lbl{mul}
\eea
and $a$ denotes the number of distinct values taken by
 $x_2, \ldots, x_k$. 

We give a coupling of $N_1$ to another random variable
$N'_1$ with the same distribution as $N_1$
 that is independent of $N_2$,
for which we can give a useful bound on
$P[N_1 \neq N'_1]$.

 Consider throwing a series of coloured balls
so each ball can land in one of the three  boxes, where the probabilities
of landing in Boxes $1,2,3$ are $1/m$, $a/m$, $(m-a-1)/m$
respectively.
First, throw $n-k$ white balls and let $N^*_1, N_2, N^*_3$ be the number of
white balls in Boxes $1,2,3$ respectively.
 Then  pick out
the balls in Boxes 1 and 3, paint them red, and throw them again. 
Then throw enough green balls so the total number of green and
red balls is $n-1$. Finally take the red balls in Box $2$ (of which
there are of $N_0$, say),
paint them blue, and throw them again but condition them to land
in Boxes 1 and 3 (or equivalently, throw each blue ball again and again until
it avoids Box 2). Then (with obvious notation, superscripts denoting
colours) set   
$$
N_1 = N_1^r + N_1^b, 
~~~ N_3 = N_3^r + N_3^b,
 ~~~ N'_1 = N_1^r + N_1^g.
$$
Then $(N_1,N_2,N_3) $ have the multinomial distribution given
by  \eq{mul}.
Also, $N'_1 $ has the same $ \Bin(n-1,\frac{1}{m})$
as $M_1$ in the statement of the lemma,
and $N'_1$ is independent of $N_2$.
Since $N'_1 = N_1 - N_1^{b} +N_1^g$,
we have that
\bean
P[ N_1 \neq N'_1] \leq
 \E[N_1^g]
+
\E[N_1^{b}]
\\
\leq
 \frac{1}{m}
 \left(k-1 + \E N_2  \right)
 + \left(
 \frac{1/m}{1 - (k-1)/m} \right) \E[N_0]
\\
\leq
\frac{1}{m}
\left( k-1 + (k-1) n/m
 +
\left( \frac{1}{1 - (k-1)/m} \right)
(k-1)n/m \right)
\eean
so that
\bean
\left|
\E[ 
(\psi_1(N_1) - \psi_1(N'_1) )  
h(N_2)
]
\right|
\leq 
P[ N_1 \neq N'_1] 
\rng(\psi_1)
 \prod_{i=2}^k \|\psi_{i}\|
\\
\leq \frac{k-1}{m} \left(
 \left( \frac{m+ n}{m} \right)
 +
\left(
\frac{
n
}{m - k+1}
\right)
\right)
\rng(\psi_1)
 \prod_{i=2}^4 \|\psi_{i}\|
\eean
and since $N'_1$ is independent of $N_2$ with
the same distribution as $M_1$, so that $\E\psi_1(N'_1) = 0$ by assumption,
$$
\E[ \psi_{1}(N'_1) h(N_2)  ]   = 0,
$$
so by \eq{0422d},
\bean
\left|
 \E \left[ W | \bX^{(k)} = \bx \right] 
 \right|
\leq \frac{k-1}{m} \left(
\left(
\frac{m+ n}{m}
\right)
 +
\left(
\frac{
n
}{m - k+1}
\right)
\right)
\rng(\psi_1)
 \prod_{i=2}^4 \|\psi_{i}\|.
\eean
Combined with \eq{0422c}, and the fact that $\|\psi_1 \| \leq \rng(\psi_1)$
since $\E[\psi_1(M_1)] =0$,
 this demonstrates the result.
$\qed$ \\

Next, we adapt Lemma \ref{momlem} to the non-uniform setting.
In this case, we need to allow $\psi_i$ to depend on the location
as well as the occupation number associated with the $i$th ball.
Consequently, some
modification of  the proof is required, and
the constants in Lemma \ref{momlem} are better than
those which would be obtained by simply 
applying the next lemma
 to the uniform case.

\begin{lemm}
\lbl{NUmomlem}
 Suppose that for $i=1,2,3,4$, $\psi_{i}$ is a real-valued function
defined on
 $ [m] \times \{0,1, \ldots,n-1\} $, 
with
 $\E[\psi_{1}(X_1,M_1) ]=0$,   set
$
\|  \psi_{i} \| : = \sup_{(x,k) \in [m ] \times [n] }
 \{ |\psi_{i}(x,k-1) |\}
$
and set
$
\rng(  \psi_{i} ) : =
 \sup_{(x,k) \in [m ] \times [n] } \{ |\psi_{i}(x,k-1) |\}
-
 \inf_{(x,k) \in [m ] \times [n] } \{ |\psi_{i}(x,k-1) |\}
$
.
Then
\bea
\left| \E \left[  \psi_{1}(X_1,M_1) \psi_{2}(X_2,M_2) \right] \right|
\leq (3 + 3 \gamma)  \rng(\psi_1) \|\psi_2\| \sum_x p_x^2 
\lbl{NU0502a}
\eea
and
\bea
\left| \E \left[ \prod_{i=1}^4 \psi_{i}(X_i,M_i) \right] \right|
\leq \left(
9 + 9 \gamma  \right)
 \rng(\psi_1) 
\left(
\prod_{i=2}^4 \|\psi_{i}\| 
\right)
 \sum_x p_x^2.
\nonumber \\
\lbl{NU0423b}
\eea
\end{lemm}


{\em Proof.} We first prove \eq{NU0423b}.
Throw $n$ balls according to the distribution $\bp$, with  four of them 
distinguished as Ball 1, Ball 2, Ball 3 and Ball 4.
For $i=2,3,4$, let $Z_i$
be the location of Ball $i$
and let $N_i$ be the number of other balls in the same urn
as Ball $i$.
Set $A = \cup_{i=2}^4 \{Z_i\}$, the union of the locations
of   Balls 2,3, and 4. 

Now suppose the balls in $A$ are painted white.
Let the balls  not in $A$  (including Ball 1 if it is not in $A$) be 
re-thrown (again, according to the distribution $\bp$).
 Those which land in $A$ when re-thrown are
painted yellow, and the others are painted red.

Now introduce
one green ball  for each white ball, and if Ball 1 is
white, let one of the green balls be labelled Ball G1.
Throw the green balls using the same distribution $\bp$.
Also, introduce  a number of blue balls equal to the number of
yellow  balls, and if Ball 1 is yellow then label
one of the blue balls as Ball B1.
 Throw the blue balls,  
but condition them to avoid $A$; that is,
use the probability mass function
$(p_x/(1 - \sum_{y \in A}p_y), x \in [m] \setminus A )$
for the blue balls.
 
Set $Z_1$ to be the location of Ball 1 (if it is  white or red)
or Ball B1 (if Ball 1 is yellow).
Set $Z'_1$ to be the location of Ball 1, if it is red or yellow, or  
the location of Ball G1 (if Ball 1 is white). 
Let $N_1^w,$ $ N_1^r$, and $ N_1^b$ respectively
denote the number of white, red, and blue
 balls at location $Z_1$, not counting Ball 1 or Ball B1 itself.
Let $N_1^y,$ $ N_1^r$, and $ N_1^g$ respectively
denote the number of yellow, red, and green
 balls at location $Z'_1$, not counting Ball 1 or Ball G1 itself.
Set
$$
N_1   = N_1^w + N_1^r + N_1^b, 
 ~~~ N'_1 = N_1^y  + N_1^r + N_1^g.
$$
Then $((Z_i,N_i)_{i=1}^4)$
 have the same joint distribution as 
 $((X_i,M_i)_{i=1}^4)$.
Also, $(Z'_1,N'_1) $ has the same distribution
as $(X_1,M_1)$ 
and $(Z'_1,N'_1)$ is independent of 
$((Z_i,N_i)_{i=2}^4)$.
Finally,  if Ball 1 is red then $Z_1 = Z'_1$ and 
 $ N'_1 = N_1 - N_1^{b} + N_1^g $,
so that
\bea
P[ (Z_1,N_1) \neq (Z'_1,N'_1)]
 \leq \E[N_1^g] + \E[N_1^{b}] + 2 P[Z'_1 \in A].
\lbl{1109a}
\eea
Now, 
\bea
P[Z'_1 \in A] \leq \sum_{i=2}^4 P [X_1 = X_i ]  = 3 \sum_x p_x^2.
\lbl{1109b}
\eea
Also, if $N_g$ denotes the number of green balls,
not including Ball G1 if Ball 1 is green,
then by \eq{gammadef},
$$
\E  [N_g]  \leq 3 + 3 n  \|\bp\| \leq 3(1 + \gamma)
$$
and also
$
 \E  [N_1^g | N_g] \leq N_g \sum_{x} p_x^2,
$  
so that
\bea
\E  [N_1^g] = \E  [ \E  [N_1^g | N_g] ] \leq 3 (1 + \gamma)     \sum_{x} p_x^2.
\lbl{1109c}
\eea
If $N_y$ denotes the number of yellow balls,
other than Ball 1,
 then 
by \eq{gammadef},
\bea
\E [ N_y ] \leq  3 n \|\bp\| \leq 3 \gamma
\lbl{1130a}
\eea
and by \eq{sixth},
\bean
 \E  [N_1^b | N_y] \leq N_y \sum_{x} 
\left( \frac{ p_x }{1 - 3 \|\bp\| } \right)^2
\leq
2 N_y \sum_x p_x^2,
\eean  
so that
\bea
\E  [N_1^b] = \E  [ \E  [N_1^b | N_y] ]  
\leq 6 \gamma 
    \sum_{x} p_x^2.
\lbl{1109d}
\eea
Set $W:= \prod_{i=2}^4 \psi_{i}(Z_i,N_i) $.
By \eq{1109a}, \eq{1109b}, \eq{1109c} and \eq{1109d},
\bea
\left|
\E[ W (\psi_1(Z_1,N_1) - \psi_1(Z'_1,N'_1) )  ]
 \right|
\leq 
P[ (Z_1,N_1) \neq (Z'_1,N'_1)] 
\rng(\psi_1)
 \prod_{i=2}^4 \|\psi_{i}\|
\nonumber \\
\leq
(9 + 9 \gamma)
 \rng(\psi_1)
 \left(
 \prod_{i=2}^4 \|\psi_{i}\|
\right)
 \sum_x p_x^2.
~~~~~~~~~
\lbl{1117a}
\eea
Since $(Z'_1,N'_1)$ is independent of 
$((Z_i,N_i)_{i=2}^4)$
 with
the same distribution as $(X_1,M_1)$, and
 $\E\psi_1(X_1,M_1) = 0$ by assumption,
$ \E[ W \psi_{1}(Z'_1,N'_1)   ]   = 0.  $
Hence,
\bean
 \E \left[ \prod_{i=1}^4 \psi_{i}(X_i,M_i) \right]
= \E \left[ W  \psi_{1}(Z_1,N_1)  \right]
= \E[ W (\psi_1(Z_1,N_1) - \psi_1(Z'_1,N'_1) )  ]
,
\eean
and then \eq{NU0423b}
follows by \eq{1117a}.
The proof of \eq{NU0502a} is similar, with the factors of 3
 replaced by 1 in \eq{1109b}, \eq{1109c} and \eq{1130a}.
$\qed$

\section{Proof of Theorems \ref{thm1} and \ref{coro1}}
\lbl{proof-1}
\allco
{\em Proof of Theorem \ref{thm1}.}
Recall the definition \eq{Yndef} of $M_i$ and $Y$.
Assume the uniform case, i.e. assume $\bp = (m^{-1},m^{-1},\ldots,m^{-1})$.
Let
$ \xi_i := {\bf 1} \{  M_i > 0 \}  $
be the indicator of the event that ball $i$ is not isolated.
Then  $Y=\sum_{i=1}^n \xi_i$, and since $\{\xi_i\}$ are exchangeable, a
random variable $Y'$ with the size-biased distribution of $Y$ can be
obtained as follows; see Lemma \ref{sblem}.
 Let $I$ be a discrete
uniform random variable over $[n]$, independent of
$X_1,\ldots,X_n$. Given the value of $I$, let 
$\bX' = (X'_1, \ldots, X'_n) \in [m]^n$
 be a random $n$-vector
with $\LL(X'_1,\ldots,X'_n) = \LL(X_1,\ldots,X_n|\xi_I=1)$.
Set
$$
Y' : = \sum_{i=1}^n {\bf 1} 
\{ \cup_{j \in [n] \setminus \{i\} }
\{X'_j= X'_i\}  \}.
$$
To apply Lemma \ref{LGthm} we need to find a random variable
$Y''$, coupled to $Y$, such that $\LL(Y'')= \LL(Y')$ and
 for some constant $B$ we have $ |Y'-Y| \le B$ (almost
surely). To check $\LL(Y'')= \LL(Y)$, we shall use  the fact that
if  $M'_I$ denotes the number of entries $X'_j$ of $\XX'$
that are equal to $X'_I$, other than $X'_I$ itself, 
then (i) given $I$ and $X_I$, $M'_I$
has the distribution of  a $\Bin(n-1,1/m)$ variable conditioned
to take a non-zero  value, and (ii) given $I, X'_I$ and $M'_I$,
the distribution of $\XX'$ is uniform over all 
possibilities consistent with the given values of $I, X'_I$
and $M'_I$. 

Define the  random $n$-vector
$
\bX := (X_1\ldots,X_n).
$
We can manufacture a random vector
$\bX'' = (X''_1,\ldots,X''_n)$,
 coupled to $\bX $ and (we assert) 
 with the same distribution
as $\bX'$,
   as follows.
\begin{itemize}
\item
Sample the random variables  $(X_1,\ldots,X_n)$.
Define $M_i$ by \eq{Yndef}.
\item
Sample a value of $I$ from the discrete uniform distribution on
 $[n]$,  independent of $\bX$.
 \item Sample a Bernoulli random variable
$\BB$ with  $P[\BB=1] = \pi_{M_I}$, where $(\pi_k, k \geq 0)$ is given by
\eq{pikdef} with $\nu=n-1$ and $p=m^{-1}$.
(By Lemma \ref{bincouplem}, $0 \leq \pi_k \leq 1$.)
 \item
Sample a value of $J$ from the discrete uniform distribution on
 $[n]\setminus \{I\}$.
 \item Define   $(X''_1,\ldots,X''_n)$ by
$$
X_i'' = \left\{ \begin{array}{ll}
X_I & {\rm if} ~ i=J ~  {\rm and} ~ \BB=1
\\
X_i & {\rm otherwise}.
\end{array}
\right.
$$
\end{itemize}

Thus $\bX''$   is obtained from $\bX$ by changing 
a
randomly selected entry of
 $\bX $ to
 the  value of  $X_I$, if $\BB=1$,
and leaving $\bX$ unchanged if $\BB=0$.

We claim that $\LL(\bX'') = \LL(\bX')$.
To see this define 
$N:= M_I$, 
and set
$N'' := -1 + \sum_{i=1}^n {\bf 1}\{X''_i= X''_I\}$.
Then $N$
has the $\Bin(n-1,m^{-1})$ distribution, while $N''$
always takes the value either $N$ or $N+1$,
taking the latter value in the case where $\BB=1$
and also $X_J \neq X_I$. Thus for any 
$k \in \{0,1,\ldots,n-1\}$, 
\bea
P[N''>k] = P[N>k] + P[N=k]  \pi_k (1- (k/(n-1))),
\lbl{1229b}
\eea 
so by the definition \eq{pikdef} of $\pi_k$, $\LL(N'')= \LL(N|N>0)$.
This also applies to the conditional distribution of $N''$
given the values of $I$ and $X_I$.

Given  the values of $N''$, $I$ and $X''_I$, the conditional
distribution of $\bX''$ is uniform over all possibilities
consistent with these given values.
Hence, 
$\LL(\bX'')=\LL(\bX')$. Therefore setting
\bean
Y'' : = \sum_{i=1}^n   {\bf 1}
 \{ \cup_{j \in [n] \setminus \{i\} } \{X''_j = X''_i\}   \},
\eean
we have that $\LL(Y'')= \LL(Y')$,
i.e. $Y''$ has the size-biased distribution of $Y$.

The  definition of $\bX''$ in terms of $\bX$ ensures that
we always have
$|Y-Y''| \leq 2$ (with equality if
    $M_I = M_J = 0$)
; this is explained further
in the course of the proof of Proposition \ref{convarprop} below.
 Thus we may apply Lemma \ref{LGthm} with $B=2$. 
Theorem \ref{thm1}
follows from that result,
 along with  the following:

\begin{prop}
\label{convarprop}
It is the case that $\Var(\E[Y''-Y|Y])  \leq
 \eta(n,m),$ where $\eta(n,m)$ is
given by \eq{etaYdef}.
\end{prop}
{\em Proof.} 
Let $\G$
be the $\sigma$-algebra generated by $\bX$.
Then $Y$ is $\G$-measurable.
By the conditional variance
formula, as in e.g. the proof of Theorem 2.1 of \cite{GP},
\bea
 \Var(\E[Y''-Y|Y]) \leq
\Var (\E[Y''-Y|\G]),
\lbl{1125a}
\eea
so it suffices to prove that
\bea
\Var(\E[Y''-Y|\G]) \leq \eta(n,m).
\label{0422a}
\eea

 For $1 \leq i \leq n$, let 
 $V_i$ denote the conditional  probability that $\BB=1$,
given $\bX$ and given that $I=i$, i.e. 
$$
V_i=\pi_{M_i}.
$$
 Let $R_{ij}$ denote the increment in the number of
non-isolated balls when the value of $X_j$ is changed to $X_i$.
Then
$$
\E[Y''-Y|\G] = \frac{1}{n(n-1)} \sum_{(i,j):i \neq j}  V_i   R_{ij}
$$
where $\sum_{(i,j):i \neq j}$ denotes summation over
pairs of distinct integers $i,j$ in $[n]$.

 For $1 \leq i \leq n$ and $j \neq i$, let
\bean
S_i := {\bf 1}\{M_i=0\}; 
~~~~~~~~~~~~~
T_i := {\bf 1}\{M_i=0\} - {\bf 1}\{M_i=1\}; \\
Q_{ij} := 
{\bf 1}\{M_i=1\}
{\bf 1}\{X_i=X_j\} .
 \eean
Then we assert that  $R_{ij}$,  the  increment in
 the number of non-isolated  balls
when ball $j$ is moved to the location of ball $i$, is given by
$
R_{ij}:= S_i + T_j + Q_{ij}.
$
Indeed, if $X_i \neq X_j$ then
 $S_i$  is the increment (if any) due to ball $i$ becoming non-isolated,
while  $T_j$ is the increment (if any)
 due either to ball $j$ becoming non-isolated,
or to another ball at the original location of ball $j$ becoming isolated when
ball $j$ 
 is moved to the location of ball $i$.
The definition of $Q_{ij}$ ensures that
if $X_i = X_j$ then   $S_i + T_j +  Q_{ij} = 0$.
 Thus,
\bea
\E[Y''-Y|\G]
 = \frac{1}{n(n-1)} \sum_{(i,j):i \neq j}  V_i(S_i+T_j +Q_{ij})
\nonumber \\
= \frac{1}{n} \sum_{i=1}^n V_i \tau_i
+ \frac{1}{n(n-1)} \sum_{(i,j):i \neq j} V_i T_j,
\lbl{0422b}
\eea
where we set
\bea
\tau_i :=
S_i + \left( \frac{1}{n-1} \right) \sum_{j :j\neq i}     Q_{ij}
= {\bf 1}\{ M_i=0\}
 + \left( \frac{1}{n-1} \right)
 {\bf 1}\{ M_i=1\}.
\lbl{taudef}
\eea
Put $a := \E[V_i]$ 
(this expectation does not
depend on $i$).  Then
by \eq{0422b}, 
\bean
\E[ Y''-Y|\G] =
\frac{1}{n} \sum_{i=1}^n \left( V_i \tau_i + aT_i   \right)
+
 \sum_{(i,j):i \neq j} 
 \frac{
(V_i-a) T_j 
}
{n(n-1)}
.
\eean
Since $(x+y)^2 \leq 2(x^2 + y^2)$ for any real $x,y$,
it follows that
\bea
\Var \left( \E[Y'-Y|\G] \right) \leq
2 \Var
\left( \frac{1}{n} \sum_{i=1}^n
 \left( V_i \tau_i  + aT_i \right)\right)
\nonumber \\
+ 2 \Var
\sum_{(i,j):i \neq j}
\frac{
 (V_i-a) T_j 
}{n(n-1)} 
.
\lbl{0513c}
\eea
From the definitions, 
 the following inequalities
hold almost surely:
\bea
- 1 \leq T_i \leq 1; ~~~ 
0 \leq V_i \leq 1; ~~~ 
0 \leq \tau_i \leq 1; 
\lbl{1229a}
\eea
and hence
\bea 
\quad \quad
\quad 
-1 \leq V_i - a  \leq 1; ~~ 
-1 \leq (V_i -a)T_j  \leq 1; ~~
- 1 \leq V_i \tau_i  + a T_i \leq 2.
\lbl{Zineqs}
\eea
Set
$ Z_i := V_i \tau_i  + a T_i, $
and $\bar{Z}_i:= Z_i - \E Z_i$. 
By \eq{Zineqs},
$
\Var Z_1 \leq \E Z_1^2 \leq 4.
$  
Also 
by \eq{Zineqs},
 we have $|\bar{Z}_i| \leq 3$,
and $-1 - \E Z_i \leq \bar{Z}_i \leq 2 - \E Z_i $.
Hence
by the case $k=2$ of Lemma \ref{momlem},
\bean
\Cov( Z_1, Z_2) = \E [ \bar{Z}_1 \bar{Z}_2]
\leq \frac{9}{m} \left(  
1 + 
\left(  \frac{n}{m-1} \right)
+  \left(  \frac{n+m}{m} \right)
\right).
\eean
Thus for the first term in the right hand side of \eq{0513c},
we have
\bea
\Var\left(  \frac{1}{n} \sum_{i=1}^n Z_i
\right)
= n^{-1} \Var(Z_1) + \left( \frac{n-1}{n}  \right) \Cov(Z_1,Z_2)
\nonumber \\
\leq  \frac{4}{n}
 +
\frac{9}{m}
 \left( 1 + \left(  \frac{n}{m-1} \right)
+  \left(  \frac{n+m}{m} \right)
\right).
\lbl{0519b}
\eea
For the second  term in the right hand side of \eq{0513c},
set $\bar{V}_i := V_i -a$. 
By \eq{1229a}, $-a \leq \bar{V}_i \leq 1-a$,
and $|T_i|\leq 1$. Hence by
the case $k=4$ of 
 Lemma \ref{momlem},
\bean
\Cov( \bar{V}_1 T_2, \bar{V}_3 T_4 )
\leq
 \E[ \bar{V}_1 T_2\bar{V}_3 T_4 ]
\leq
 \frac{3}{m} \left( 2 + \left( \frac{n}{m-3}
\right) + 
 \frac{n}{m}
 \right).
\eean
By 
\eq{Zineqs}, we can always bound
 $\Cov(\bar{V}_i T_j,\bar{V}_{i'} T_{j'})$
by $1$. Hence, expanding $\Var \sum_{(i,j): i \neq j} \bar{V}_i T_j$
in the same manner as with \eq{0519a} below, 
yields
\bean
\quad 
  \Var
 \sum_{(i,j): i \neq j}
 \frac{
 \bar{V}_i  T_j 
}{n(n-1)}
\leq 
 \frac{3}{m} \left( 2 + \left( \frac{n}{m-3}
\right) +
 \frac{n}{m}
 \right)
+ \frac{4}{n} +  \frac{2 }{n(n-1) } .
\eean
Using this with \eq{0513c} and \eq{0519b}
yields
\bean
 \Var \left( \E[Y'-Y|\G] \right)
\leq
\frac{16}{n} + \frac{4}{n(n-1)} + \frac{24}{m} 
\left( 2 + \frac{ n}{m-3} + \frac{n}{m} \right) . 
\eean
 This completes the proof of
Proposition \ref{convarprop}, and hence of Theorem \ref{thm1}.
 $\qed$ \\

{\em Proof of Theorem \ref{coro1}.}
Suppose $n,m$ both go to infinity in a linked manner, in
such a way that $n/m \to \alpha \in (0,\infty)$.
Then
 it can be shown (see Theorem II.1.1 on pages 37-38 of
\cite{KSS}) that
$
\E Y \sim n (1 - e^{-\alpha}) ,
$
 and
\bean
\Var(Y)
\sim n \left(  e^{-\alpha} (1- e^{-\alpha}) 
+ 
 e^{-2 \alpha}(
 \alpha(1 - \alpha))
\right) 
= n g(\alpha)^2.
 \eean
Substituting 
these asymptotic expressions
 into \eq{thm1eq} and using the fact that
in this asymptotic regime,
$
\left( n \eta(n,m) \right) \to
 16 + 24 \alpha ( 2 + 2 \alpha  ), 
$
yields \eq{0901d}.

\section{The non-uniform case: proof of Theorem \ref{varthm}}
\lbl{secvarpf}
\allco
For this proof,  we use the following
notation.  Given $n$, $m$,  and the probability
 distribution $\bp$ on $[m]$, 
let $X_1,X_2,\ldots,X_{n+1}$ be independent $[m]$-valued random variables
with common probability mass function $\bp$.
Given  $ i \leq j \leq n +1$,  set
$\bX_j:= (X_1,\ldots,X_j)$ and
$$
\XX_{j \setminus i} := 
\left\{
\begin{array}{lll}
(X_1, X_2, \ldots, X_{i-1}, X_{i+1}, \ldots, X_j)
& {\rm if} & 1 < i < j 
\\ ( X_2, \ldots,  X_j) & {\rm if} &  i = 1 \\
 ( X_1, \ldots,  X_{j-1}) & {\rm if} &  i = j.
\end{array}
\right. 
$$
Given any sequence $\bx = (x_1,\ldots, x_k)$, set
\bea
H(\bx) =\sum_{i=1}^k \left( 1 - \prod_{j\in [k] \setminus \{i\} }
 (1 - {\bf 1}  \{ x_j = x_i   \} ) \right),
\lbl{Hdef}
\eea
which is the number of non-isolated entries in the sequence $\bx$, so that
in particular, $Y= H(\XX_n)$.
We shall use the following 
 consequence of  Jensen's inequality:
for all $k \in \N$,
\bea
\lbl{Jensen}
\quad
\quad
\quad
(t_1 +t_2 + \cdots + t_k)^2 \leq k (t_1^2 + \cdots + t^2_k ), 
\quad
 \forall 
 ~~(t_1,\ldots,t_k) \in \R^k. 
\eea
We shall also use several times the fact that $- t^{-1} \ln (1-t)$
is increasing on $t \in (0,1)$ so that by \eq{sixth}
for all
$x \in [m]$ 
 we have 
\bea
\ln (1-p_x ) \geq ( 11 \ln (10/11) ) p_x \geq  - 1.05 p_x
\lbl{qlb}
\eea
whereas $(1-e^{-t})/t$ is decreasing on $t \in (0,\infty)$ so that
by \eq{gammadef},
for any $\alpha >0$  and $x \in [m]$ we have
\bea
1 - e^{-\alpha n p_x} \geq
 (1- e^{-\alpha \gamma})
(n p_x / \gamma).
\lbl{qub}
\eea

{\em Proof of \eq{genvarub}.}
We use Steele's variant of the Efron-Stein inequality \cite{StES}.
This says, among other things, that when (as here) $X_1,\ldots,X_{n+1}$
are independent and identically distributed random variables
and $H$ is a symmetric function on $\R^n$,
\bean
\Var H (X_1,\ldots,X_n) \leq 
\frac{1}{2} \sum_{i=1}^n
\E  [ (H( \XX_n)  - H(\XX_{(n+1) \setminus i} ))^2 ] 
\\
 =
(n/2)
\E [  ( H( \XX_n ) - H(\XX_{(n+1) \setminus n}) )^2  ].
\eean
Hence, by the case $k=2$ of \eq{Jensen},
\bean
\Var Y \leq
 n  ( \E [ (H(\bX_n) - H(\bX_{n-1}))^2 ] +  
 \E [ (H(\bX_{(n+1) \setminus n}) - H(\bX_{n-1}))^2 ] 
)
\\  
 =2 n   \E [ (H(\bX_n) - H(\bX_{n-1}))^2 ]. 
\eean  
With $M_j$ defined by \eq{Yndef}, $H(\bX_n) - H(\bX_{n-1})$ 
is equal to ${\bf 1}\{M_n \geq 1\} + 
{\bf 1}\{M_n = 1\} $, so is nonnegative and bounded by 
 $2 {\bf 1}\{M_n \geq 1\} $. 
Therefore,
$$
\Var [Y] \leq 8 n P [M_n \geq 1] \leq 8 n \E M_n \leq 8n^2 \sum_x p_x^2.
~~~
\qed
$$  
%

{\em Proof of  \eq{genvarlb}.}
Construct a martingale as follows.
Let $\F_0 $ be the trivial $\sigma$-algebra, and
for $i \in [n]$ 
let $\F_i := \sigma(X_{1},\ldots,X_{i})$ and write $\E_i$ for
conditional expectation given $\F_i$.
Define martingale differences
$
\Delta_i = \E_{i+1} Y - \E_{i} Y.
$
Then $Y - \E Y = \sum_{i=0}^{n-1} \Delta_i$,
and by orthogonality of martingale differences,
\bea
\Var [Y] = \sum_{i=0}^{n-1} E[ \Delta_i^2] =
  \sum_{i=0}^{n-1} \Var [ \Delta_i].
\label{mardif}
\eea
We look for lower bounds for $E[\Delta_i^2]$. 
Note that
\bea
  \Delta_i = \E_{i+1}[ W_i] , 
~~~~
{\rm where} 
~~~~
W_i : = H(\XX_{n}) - H(\XX_{(n+1) \setminus {(i+1)}}).
\lbl{1114a}
\eea 
Recall from \eq{Yndef} that for $i < n$,
 $M_{i+1}$ denotes the number
 of balls in the sequence of $n $ balls,  other than ball $i+1$, in
the same position as ball $i+1$. Similarly,
define $M_{n+1}$ and 
$M_k^i$  
(for $k \in [n+1]$)
by
\bea
\quad \quad
M_{n+1} : = \sum_{j \in [n] }  {\bf 1} \{ X_j = X_{n+1} \} ; 
~~~~~~
M^i_{k} : = \sum_{j \in [i] \setminus \{k\} }  {\bf 1} \{ X_j = X_{k} \} , 
\lbl{MMdef}
\eea
so that $M_{n+1}$ is the number of
balls, in the sequence of $n$ balls, in the same location as ball $n+1$,
while $M^i_k$ is similar to $M_k$,
but defined in terms of
 the first $i$ balls, not the first $n$ balls.

Set $h_0(k) := {\bf 1}\{k\geq 1\} + {\bf 1}\{k = 1\}$.  Then
$ H ( \bX_n ) - H(\bX_{n\setminus (i+1)} ) =  h_0(M_{i+1}) , $
and if $X_{n+1} \neq X_{i+1}$ then
$
H ( \bX_{(n+1) \setminus (i+1)} ) - H(\bX_{n\setminus (i+1)} )
 =  h_0(M_{n+1}) , 
$
so that $W_i = h_0(M_{i+1}) - h_0(M_{n+1})$ in this case.
For taking $\E_{i+1}$-conditional expectations,
it is convenient to approximate 
 $ h_0(M_{i+1}) $ and  $h_0(M_{n+1})$ by 
 $ h_0(M^i_{i+1}) $ and  $h_0(M_{n+1}^i)$ respectively.
To this end, define 
\bea
Z_i : = W_i -  
( h_0(M^i_{i+1}) - h_0(M^i_{n+1}) ).
\lbl{Zdef}
\eea
Since $ h_0(M^i_{i+1}) $ 
is $\F_{i+1}$-measurable,
taking 
 conditional expectations
yields
\bea
   h_0(M^i_{i+1})  
=
\E_{i+1} [ W_i ] 
+ \E_{i+1} [ h_0 (M^i_{n+1} ) ]
- \E_{i+1} [Z_i].
\lbl{1109g}
\eea 
Set $ \delta : = (288 \gamma e^{1.05 \gamma} )^{-1}.  $
 We shall show that for $i$ close to $n$, 
in the sense that $n-\delta n \leq i \leq n$,
the variances of the terms on the right
of \eq{1109g}, other than
 $\E_{i+1} [W_i]$,
are small compared to the variance of the left hand side,
 essentially because
$\E_{i+1} [ h_0(M_{n+1}^i) ]$
 is more smoothed out than 
$h_0(M_{i+1}^i) $, while 
$P[Z_i \neq 0]$
 is small when  $i$ is close to $n$.
These estimates then 
yield a lower bound on the variance of 
 $\E_{i+1} [W_i]$.

First consider the left hand side $h_0(M^i_{i+1}) $.
This variable takes the value 0 when
$M_{i+1}^i =0 $, and takes a value at least $1$
when $M_{i+1}^i \geq 1 $.  Hence, 
\bea
 \Var [h_0(M_{i+1}^i) ] \geq 
(1/2)
 \min ( P[M_{i+1}^i = 0 ],
P[M_{i+1}^i \geq 1 ]).
\label{1109e}
\eea
For $i \leq n$, by \eq{qlb} and \eq{gammadef},
 \bea
P[M_{i+1}^i = 0] =
\sum_x p_x(1-p_x)^{i} 
\geq
\sum_x p_x (1-p_x)^{n} 
 \nonumber \\
\geq  \sum_x     p_x e^{-1.05 n p_x} 
\geq \gamma^{-1}
e^{-1.05 \gamma}
  \sum_x n p_x^2.
\lbl{1117b}
\eea
For $i \geq (1- \delta)n$ we have $i  \geq n/2$, 
so by \eq{qub} and the fact that $\gamma \geq 1$ by  \eq{gammadef},
 \bea
P[M_{i+1}^i \geq 1] =
\sum_x p_x( 1 - (1-p_x)^{i} )
\geq
\sum_x p_x( 1 -  e^{-np_x/2})
\nonumber \\
\geq
\sum_x p_x( 1 -  e^{- \gamma /2}) np_x/\gamma
\geq (1-e^{-1/2})\gamma^{-1} \sum_x   n  p_x^2.
\lbl{1117c}
\eea
Since $\gamma \geq 1$, and $e^{-1.05} < 1 - e^{-0.5}$,
the lower bound in \eq{1117b} is always less than that in \eq{1117c},
so combining these two
estimates and using \eq{1109e} yields 
\bea
\quad \quad
 \Var [h_0(M^i_{i+1}) ] \geq (1/2) \gamma^{-1} e^{-1.05 \gamma} 
 \sum_x n p_x^2, ~~~~~
i \in [n-\delta n,n].
\lbl{1109f}
\eea
Now consider the second term $\E_{i+1} [ h_0(M_{n+1}^i) ]$
 in the right hand side of \eq{1109g}.
Set $N_x^i := \sum_{j=1}^i {\bf 1} \{ X_j =x \}$, and
for $1 \leq \ell   \leq i$ set
$M_\ell^i$ to be $N_{X_\ell}^i -1$. Also
 set
$\tih_0(k) = (k+1)^{-1}  h_0 (k+1)$.
Then, since $h_0(0)=0$, we have that
\bea
      \Var ~   \E_{i+1} [h_0(M_{n+1}^i)]  
 = \Var \sum_x p_x  h_0( N_x^i) 
 = \Var \sum_{j=1}^i p_{X_j} \tih_0( M_j^i) 
\nonumber  \\ 
       =  
\frac{i}{n^2} \Var [ np_{X_1} \tih_0(M_1^i) ]
+  \frac{i(i-1)}{n^2} 
\Cov \left[ np_{X_1} \tih_0(M_1^i), np_{X_2} \tih_0(M_2^i) \right].
\lbl{1117d}
\eea
Suppose $i \leq n$.
Since $0 \leq n p_{X_1} \tih_0 (M_1^i)  \leq 2 np_{X_1}$,
\eq{gammadef} yields
\bean  
\frac{i}{n^2} \Var [ np_{X_1} \tih_0(M_1^i) ]
\leq
n^{-1} \E [ 4n^2 p_{X_1}^2 ] = 4n \sum_x p_x^3 
\leq 4 \gamma \sum_x p_x^2,
\eean
while by Lemma \ref{NUmomlem} and \eq{gammadef}, 
\bean
\Cov \left[ np_{X_1} \tih_0(M_1^i), np_{X_2} \tih_0(M_2^i) \right]
\leq
(3 + 3 \gamma) 4 \gamma^2 \sum_x p_x^2.
\eean
Combining the last two estimates on \eq{1117d} 
and using assumption \eq{nbig1} yields
\bea
\Var  ( \E_i [h_0(M_{n+1}^i)] ) \leq  
( 1 + 3 \gamma + 3  \gamma^2)4 \gamma \sum_x p_x^2 
\nonumber \\
\leq
(18 \gamma e^{1.05 \gamma})^{-1}
n \sum_x  p_x^2.
\lbl{1109i}
\eea 
We turn to the third term in the right hand side of \eq{1109g}.
As discussed just before 
 \eq{Zdef},  when $X_{n+1} \neq X_{i+1}$ we have
$W_i = h_0(M_{i+1} ) - h_0(M_{n+1})$, and
it is clear from the definitions \eq{1114a} and 
\eq{MMdef}
that if  $X_{n+1} = X_{i+1}$  
then both $W_i$ and $
h_0(M_{i+1}^i) - h_0(M_{n+1}^i) $ are zero, and therefore by \eq{Zdef}, 
$$
Z_i = 
(
  h_0(M_{i+1} ) - h_0(M_{i+1}^i)
 - h_0(M_{n+1})
 + h_0(M_{n+1}^i) 
 ) {\bf 1}\{X_{n+1} \neq X_{i+1} \}.
$$
By the conditional Jensen inequality,
\bean
\Var ( \E_{i+1} [ Z_i  ] ) 
\leq \E [ ( \E_{i+1} [ Z_i ] )^2 ] 
\leq \E [ Z_i^2 ].
\eean
The random variable
 $  h_0(M_{n+1} ) - h_0(M_{n+1}^i)  $ lies in the range $[-2,2]$
 and  is zero unless 
$X_j = X_{n+1}$ for some $j \in  (i,n]$.
Similarly,
$  h_0(M_{i+1} ) - h_0(M_{i+1}^i)  $ lies
in $[-2,2]$ and is zero unless
$X_j = X_{i+1}$ for some $j \in  (i+1,n]$.
Hence, using \eq{Jensen} and the definition of $\delta$
yields for $i \in [n-\delta n,n]$ that
\bea
\Var ( \E_{i+1} [ Z_i  ] ) 
\leq 2 ( 4 P[M_{n+1} \neq M_{n+1}^i] + 4 P[M^i_{i+1} \neq M^i_{i+1}]) 
 \nonumber \\
\leq 16 \delta n \sum_x p_x^2
\leq (18 \gamma e^{1.05 \gamma} )^{-1}  \sum_x n p_x^2.
\label{1109h}
 \eea

By \eq{1109g} and the case $k=3$ of \eq{Jensen}, 
\bean
\Var [  h_0 (M_{i+1}^i) ]  \leq  3 ( \Var (\E_{i+1} [W_i]) +
\Var ( \E_{i+1}[ h_0(M_{n+1}^i) ] ) 
+
 \Var ( \E_{i+1}[
Z_i ] ) 
).
\eean
Rearranging this and using \eq{1109f}, \eq{1109i}, and \eq{1109h} yields
 the lower bound
\bean
 \Var (\E_{i+1} [W_i])
\geq \left(
\frac{1}{6} - \frac{2}{18}
 \right)  \frac{ e^{-1.05 \gamma} }{\gamma} \sum_x n p_x^2
= 
 \frac{e^{-1.05 \gamma}}{18 \gamma}
 \sum_x n p_x^2, 
\eean
for $ i \in [n-\delta n, n]$.
Since the definition of $\delta$, the condition
\eq{nbig1} on $n$ and the assumption $\gamma \geq 1 $ guarantee
that $n \delta \geq 2$,  and since $\lfloor t \rfloor \geq 2t/3$ for 
$t  \geq 2$, 
 by \eq{mardif} and \eq{1114a} we have
$$
\Var [Y] \geq 
\lfloor \delta n \rfloor (18 \gamma e^{1.05 \gamma})^{-1} n \sum_x p_x^2
\geq (\delta n ) (27 \gamma e^{1.05 \gamma})^{-1} n \sum_x p_x^2
$$
which is \eq{genvarlb}. $\qed$

\section{Proof of Theorem \ref{nonunithm}}
\allco 
{\em Proof of \eq{thm2eq2}.}
Write $\sigma $ for $\sigma_Y$, and for $t \in \R$
set $F(t):= P[(Y-\E Y)/\sigma \leq t]$.
Set $z_0 := \sigma^{-1} ( \lfloor \E Y \rfloor - \E Y)$, and set $z_1 := z_0 + 
(1 - \eps )/\sigma$, for some $\eps \in (0,1)$.  Then
since $Y$ is integer-valued, 
$F (z_1) = F (z_0) $.
On the other hand, by the unimodality of the normal density,
\bean
\Phi(z_1) - \Phi (z_0) 
\geq
(1 - \eps) \sigma^{-1} (2 \pi)^{-1/2} \exp( - 1/ (2 \sigma^2))
\eean
so that
$ D_Y  $
is at least half the expression above. Making $\eps \downarrow 0$
and using the fact that $e^{-1/(2 \sigma^2)} \geq e^{-1/2}$ for
$\sigma \geq 1$, 
gives us \eq{thm2eq2} in the case where $\sigma \geq 1$.

When $\sigma < 1$, we can take $ z_2 \leq 0 \leq z_3$,
with $z_3 = z_2 +1$ and
$F (z_3)  = F (z_2) $.
By the $68-95-99.7$ rule for the normal distribution, 
$\Phi (z_3) - \Phi (z_2) \geq 1/3 $, so 
$ D_Y  \geq 1/6$, giving
 us \eq{thm2eq2} 
in the case where  $\sigma < 1$.
 $\qed$ \\

%

 So in Theorem \ref{nonunithm}, the difficulty lies entirely in proving
the upper bound in \eq{thm2eq}, under assumptions \eq{sixth} and \eq{nbig1}
which we assume to be in force throughout the sequel.
By \eq{nbig1} we always have $n \geq 1661$.

As before, set $h_0(k): = {\bf 1}\{k \geq 1\} + {\bf 1}\{k=1\}$.
Define for nonnegative integer $k$ the functions
\bea
h_1(k) := 
1 - h_0(k) =
 {\bf 1}\{k= 0\} - {\bf 1}\{k= 1\};
~~
~~
\nonumber
\\
h_2(k) : = 2 {\bf 1}\{k= 1\}
- {\bf 1}\{k= 2\}; 
~~ ~~
~~ ~~
h_3(k):= {\bf 1}\{k=1\}.
\nonumber
\eea
The function $h_1(k)$
may be interpreted as the increment in 
the number of non-isolated balls should a ball in an  urn containing
$k$ other balls be removed from that urn with the removed
ball then deemed to be non-isolated itself.  If $q=0$ then the ball
removed becomes non-isolated so the increment is $1$,
while if $q=1$ then the other ball in the urn becomes
isolated so the increment is $-1$.

The function $h_2(k)$ is chosen so that 
$h_2(k) + 2 h_1(k)$
(for $k \geq 1$)
is the increment in the number of non-isolated balls
if two balls should be removed from an urn containing $k-1$ other balls,
with both removed balls  deemed to be non-isolated.
The interpretation of $h_3$ is given later.

We shall need some further functions $h_i$ which we
define here to avoid disrupting the argument later
on.
For $x \in [m]$ and $k \in \{0\} \cup [n-1]$, let $\pi_k(x)$ be 
given by the $\pi_k$ of \eq{pikdef} when $\nu = n-1$ and $p=p_x$. 
With the convention $0 \cdot \pi_{-1}(x) := 0 \cdot
 h_i(-1) : = 0$,  define
\bean
h_4(k,x) : = \frac{k \pi_{k-1}(x)}{n-1} 
+ \frac{(n-k-1) \pi_{k}(x)}{n-1} 
- 1 ; \\
h_5(k,x) := \pi_k(x)/(n-1), ~~~ 
h_6(k) := k h_2(k)  ; ~~~ 
\\
h_7(k,x) : = h_3(k)
+ h_4(k,x)  
 - \frac{k (2 + h_4(k,x) ) h_1(k-1)}{n}  \\
- \frac{k h_5(k,x)
 (k-1) h_2
 (k-1)}{n}.
 \eean
For $i = 0,1,2,3,6$ define $h_i(k,x):= h_i(k)$. 
For each $i$ define
 \bea
\tih_i(k,x):= h_i(k+1,x)/(k+1). 
\lbl{tihdef}
\eea
Sometimes we shall write
 $\tih_i(k)$ for $\tih_i(k,x)$ when $i \in \{0,1,2,3,6\}$.
Define  $\|h_i\| := \sup_{k,x} |h_i(k,x)|$ and
  $\|\tih_i\| := \sup_{k,x} |\tih_i(k,x)|$.

Now we estimate
some of the $h_i$ functions.  
Since 
$\pi_0(x) = 1$ we have $h_4(0,x) = h_7(0,x) =0$ for all $x$,
which we use later.
Also,
by Lemma \ref{bincouplem},
\bea
- 1 \leq h_4(k,x) \leq 0
\lbl{1123a}
\eea 
and  $h_7(1,x) = 1 + h_4(1,x) (1 -n^{-1}) - 2/n$ so that
$ -1/n \leq h_7(1,x) \leq 1. $
Also, since \eq{nbig1} implies $n \geq 1661$,
$$
h_7(2,x) 
= h_4(2,x)(1+2 n^{-1}) - \frac{4 \pi_2(x)}{n(n-1)} + \frac{4}{n} 
\in [-1,4/n].
$$
Also, $h_3(3) = h_1(2) =0$ so that by \eq{1123a}
\bean
h_7(3,x) 
 =  \frac{6 \pi_3(x)  }{n(n-1) } + h_4(3,x)
\in [-1,1],
\eean
again since $n \geq 1661$.
For $k \geq 4$, $h_7(k,x) = h_4(k,x) \in [-1,0]$. Thus,
\bea
\| h_7 \| \leq 1;
 ~~~
\| \tih_7 \| \leq
 1; ~~~
\lbl{1125b1}
\\
 \|\tih_0\| =2; ~~~
\| h_5 \| \leq (n-1)^{-1}; ~~~
\| h_6 \| = 2.
\lbl{1125b2}
\eea

The strategy to prove Theorem \ref{nonunithm} is similar to
 the one used in the uniform case, but the construction of 
a random variable with the distribution of $Y'$, where $Y'$ is defined
 to have the $Y$ size biased distribution,
is more complicated. As in the earlier case,
by Lemma \ref{sblem},
if $I$ is uniform over $[n]$ then
the distribution of the sum  $Y$ conditional on 
$M_I>0$ has the  distribution of $Y'$. 
However, in the non-uniform case the
conditional information that $M_I > 0$
affects the distribution of $X_I$. Indeed,
for each $i$, by Bayes' theorem
\bea
P[X_i = x | M_i > 0] = \frac{p_x (1-(1-p_x)^{n-1})}{\sum_y p_y (1-(1-p_y)^{n-1})}
=: \phat_x.
\lbl{phatdef}
\eea
Therefore the conditional distribution of  $(X_1,\ldots,X_n)$, 
given that $M_i > 0$, is  obtained by sampling $X_i$ with
probability mass function $\phat$ and then sampling
$\{X_j,j \in [n] \setminus \{i\}\}$ independently with probability mass
function $p$, conditional on at least one of them taking
the same value as $X_i$.  Equivalently, 
sample the value of $X_i$, then  $M_i$ 
according
to the  binomial $\Bin(n-1,p_{X_i})$
 distribution conditioned to be at least one,
then select a subset ${\cal J}$ of $[n]\setminus \{i\}$ uniformly at random 
from sets of size $M_i$,  let the values of $X_j, j \in {\cal J}$
be equal to $X_i$, and let
 the values of $X_j, j \notin {\cal J}$
be independently sampled from the distribution
with the probability mass function of $X$ given that
$X \neq X_i$.

Thus a random variable $Y''$, coupled to $Y$ and
with the same distribution as $Y'$,
 can be obtained as follows.
First sample $X_1,\ldots,X_n$ independently from the original
distribution $\bp$, and set $\bX = (X_1,\ldots,X_n)$; then select $I$ uniformly at random from $[n]$.
Then sample a further random variable $X_0$  with the 
probability mass function $\phat$. Next, change the value of
$X_I$ to that of $X_0$; next let $N$ denote the number of other
values $X_j, j \in [n] \setminus \{ I \}$ which are equal to $X_0$,
and let $\pi_k = \pi_k(X_0)$ be defined by \eq{pikdef} with
$\nu = n-1$ and $p = p_{X_0}$.  Next, sample
a Bernoulli random variable ${\cal B}$ with
parameter $\pi_N$, and
if ${\cal B}=1$   change the value
of one of the $X_j, j \in [n] \setminus \{I\}$
 ($j=J$, with $J$ sampled uniformly
at random from all possibilities) to $X_0$.    
Finally, having made these changes, define $Y''$ in
the same manner as $Y$ in  the original sum \eq{Yndef} but in
terms  of the changed variables.
Then $Y''$ has the same distribution as $Y'$ by
a similar argument to that given around \eq{1229b}
in the uniform case.

Having defined coupled variables $Y,Y''$ such that
$Y''$ has the $Y$ size biased distribution, we wish
to use Lemma \ref{LGthm}. To this end, we need to estimate
the quantities denoted $B$ and $\Delta$ in that lemma.
The following lemma makes a start.
Let $\G$ be the $\sigma$-algebra generated by the value
of $\XX$, and for $x \in [m]$ 
let $N_x:= \sum_{i=1}^n {\bf 1}\{X_i = x\}$ be the number of  balls in urn $x$.
\begin{lemm}
\lbl{1125lem}
It is the case that
\bea
|Y'' -Y | \leq 3, ~~~a.s.
\lbl{1114b}
\eea
and 
\bea
\E[Y''-Y |\G] =   
 2 + \left( \sum_x \hat{p}_x
 h_5(N_x,x)  n^{-1} \sum_{i=1}^n h_6(M_i) 
\right)
+ 
 \left(  \sum_x \phat_{x}  h_7(N_x,x)  \right)
\nonumber \\
- \left(  \sum_x \phat_{x}  h_{4}(N_x,x)  \right)
\left( n^{-1} \sum_{i=1}^n h_0(M_i) \right)
- \frac{2}{n} 
 \sum_{i=1}^n h_0(M_i) 
. ~~~~~~~~~~
\label{1118c}
\eea
\end{lemm}
{\em Proof.}
We have
\bea
\E[Y''-Y |\G] = 
\E[Y''-Y |\XX] = 
\sum_x  \phat_x
\E_x[Y''-Y |\XX],  
\lbl{CE1}
\eea
where $\E_x [\cdot|\XX]$ is conditional expectation
given the value of $\XX$ and given also that $X_0=x$.
%
The formula for
$\E_x[Y''-Y|\XX]$ will depend on  $x$  through the value of $N_x$
and through the value of $p_x$.

We distinguish
between  the cases where $I$ is selected with $X_I=x$ (Case I) and
where $I$ is selected with $X_I \neq x$ (Case II).
If $N_x=k $,
then in Case I the value of $N$ on which is based
the probability $\pi_N(x)$ of importing a further ball to $x$
is $k-1$ whereas in  Case II this value of $N$ is $k$.
The probability of Case I occurring is $k/n$.

The increment $Y''-Y$   gets a  contribution 
of $h_1(M_i)$ from the  moving of Ball
$i$ to $x$ in Case II, and gets a further contribution
of $h_1(M_j) + h_2(M_i){\bf 1}\{ X_i= X_j\}$ if  $X_j $ is also
imported to $x$ from a location distinct from $x$.
Finally, if $N_x =k$  the increment gets a further contribution of 
$h_3(k)$
from the fact that if there is originally a single ball  at $x$,
then this ball will no longer be isolated after importing at
least one of balls  $I$ and $J$ to $x$ 
(note that $\pi_0(x) = 1$ so we never end up with an isolated
ball at $x$).
Combining these contributions, we have  \eq{1114b}, and also that
 for values of $\XX, x$ with $N_x=k$,
\bean
\E_x[Y''-Y|\XX] =
 h_3(k) + \frac{k \pi_{k-1}(x)}{n(n-1)}
\sum_{\{j:X_j \neq x\}} h_1(M_j)
+  n^{-1}
\sum_{\{i:X_i \neq x\}} h_1(M_i)
\nonumber \\
+
  \frac{ \pi_{k}(x)}{n(n-1)}
  \sum_{(i,j):i \neq j, X_i \neq x, X_j \neq x }
( h_1(M_j) +
h_2(M_i) {\bf 1}\{X_i= X_j\} ) 
\eean
where
in the right hand side, the first sum comes from Case I
and the other two sums come from Case II.
Hence, if $N_x=k$ then
\bea
\E_x[Y''-Y|\XX] =
 h_3(k) 
+ \left( \frac{\pi_k(x)}{n(n-1)}
 \sum_{\{i:X_i \neq x\}} M_i h_2(M_i) \right)
\nonumber \\
+ \left(
\frac{k \pi_{k-1}(x)}{n(n-1)} + \frac{1}{n} 
+ \frac{(n-k-1) \pi_{k}(x)}{n(n-1)}
\right)
\sum_{\{j:X_j \neq x\}} h_1(M_j)
\nonumber \\
= h_3(k) + \frac{2 + h_4(k,x)  }{n} 
\left(
\left(
 \sum_{i=1}^n h_1(M_i)
\right)
 - k h_1(k-1) 
\right)
\nonumber \\
+ \frac{ h_5(k,x) }{n} 
\left(
\left(
 \sum_{i=1}^n h_6(M_i)
\right)
 - k h_6(k-1) 
\right)
\nonumber \\
 = 2 + h_7(k,x) - 
\left( \frac{h_{4}(k,x)}{n}
 \sum_{i=1}^n h_0(M_i) 
\right)
- \left( \frac{2}{n}
 \sum_{i=1}^n h_0(M_i) \right)
\nonumber  \\
+
 \frac{h_5(k,x)}{n} \sum_{i=1}^n h_6(M_i) 
.
\nonumber
\eea
Then by 
\eq{CE1}
 we have \eq{1118c}. $\qed$ \\

The next lemma is based on the observation that since 
$h_4(0,x) = h_7(0,x)=0$, the sums 
of the form $\sum_x$ in \eq{1118c} are over non-empty urns
so can be expressed as sums over the balls, i.e. of the
form $\sum_{i=1}^n$. 
We need further notation.
Set $\xi_i :=  \phat_{X_i} \tih_{4}(M_i,X_i)$, and 
  $T_j:= h_0(M_j)$. Set
 $b := \E[T_j]$ (this 
does not
depend on $j$), and
 $\bar{T}_j:= T_j - b$.
Again write  $\sum_{(i,j):i \neq j}$
for $\sum_{i=1}^n \sum_{j \in [n] \setminus \{i\}}$.
\begin{lemm}
It is the case that
\bea
\Var (\E[Y''-Y |Y]) \leq 
  12 (n-1)^{-2}
+ 3 n^{-2} \Var \left(
 \sum_{(i,j): i \neq j} 
\xi_i \bar{T}_j  \right)
\nonumber \\
+  3   \var   \sum_{i=1}^n
 (  
  [\tih_7(M_i,X_i) - (1-n^{-1}) b \tih_4(M_i,X_i)  ]
 \phat_{X_i}
\nonumber \\
-
[ n^{-1} 
 h_0(M_i) ( 2 + \phat_{X_i} \tih_{4}(M_i,X_i) )
]
)
.
\lbl{1118b}
\eea
\end{lemm}
{\em Proof.}
As in Section \ref{proof-1},
\eq{1125a} holds here too.
 So it suffices to prove \eq{1118b} with the left hand
 side replaced by $\Var(\E[Y''-Y |\G])$.
Set 
\bea
\rho(\bX) :=
 2 + \sum_x \hat{p}_x
h_5(N_x,x)
n^{-1} \sum_{i=1}^n h_6(M_i) .
\lbl{rhodef}
\eea

Using \eq{tihdef},
we reformulate the sums in \eq{1118c} as follows.
Since
$h_7(0,x) = 0$,
\bea
   \sum_x \phat_{x}  h_7(N_x,x)  
=
 \sum_{i=1}^n  \phat_{X_i} \left( \frac{  h_7(M_i +1,X_i )  }{ (M_i+1) } 
\right)
=
 \sum_{i=1}^n
 \phat_{X_i}  \tih_7(M_i,X_i ).   
\lbl{1118g}
\eea   
Similarly,
 $  \sum_x \phat_{x}  h_{4}(N_x,x)  =
 \sum_{i=1}^n \phat_{X_i}  \tih_{4}(M_i,X_i ) $ so that  
\bea
 \left(  \sum_x \phat_{x}  h_{4}(N_x,x)  \right)
 n^{-1} \sum_{i=1}^n h_0(M_i) 
= \left(n^{-1} \sum_{i=1}^n
 \phat_{X_i} \tih_{4}(M_i,X_i) h_0(M_i) 
\right)
\nonumber \\
+
n^{-1} \sum_{(i,j):i \neq j} \phat_{X_i} \tih_{4}(M_i,X_i) h_0(M_j). 
~~~~
\label{0931c}
\eea
Substituting \eq{rhodef}, \eq{1118g} and 
\eq{0931c}  into \eq{1118c} gives
\bea
\E[Y''-Y |\G] =   
 \left(  \sum_{i=1}^n \phat_{X_i}  \tih_7(M_i,X_i)  \right)
- \left( 
 \sum_{i=1}^n
\frac{ h_0(M_i) }{n} ( 2 + \phat_{X_i} \tih_{4}(M_i,X_i) )
\right)
\nonumber \\
+ 
\rho(\bX)
- n^{-1} \sum_{(i,j):i \neq j} \phat_{X_i} \tih_{4}(M_i,X_i) h_0(M_j). 
      ~~~~~~~~~~~~~~~~~~~~~~~
\lbl{1125f}
\eea
 The last sum in \eq{1125f} can be rewritten as follows:
\bean
\sum_{(i,j):i \neq j}
 \xi_i T_j 
=
\sum_{(i,j):i \neq j}
 \left( b \xi_i  
+ \xi_i( T_j -b)
\right)
= 
(n-1)
\left( 
 \sum_{i=1}^n b \xi_i
\right)
+
 \sum_{(i,j):i \neq j} 
\xi_i  \bar{T}_j  .
\eean
Substituting into \eq{1125f} yields
\bean
\E[Y''-Y |\G] =   
\rho(\bX)
- 
 \left(
n^{-1}
 \sum_{(i,j): i \neq j} 
\xi_i \bar{T}_j \right)
+
\sum_{i=1}^n
 (  
 \phat_{X_i}  \tih_7(M_i,X_i)
\nonumber \\
- 
(1-n^{-1}) b \phat_{X_i} \tih_4(M_i,X_i)  
-
[ n^{-1} 
 h_0(M_i) ( 2 + \phat_{X_i} \tih_{4}(M_i,X_i) )
]
).
\eean
By \eq{rhodef} and \eq{1125b2},
 $|\rho(\bX) -2| \leq 2(n-1)^{-1} $,
so that $\Var ( \rho(\bX)) \leq 4(n-1)^{-2} $.
By 
 \eq{Jensen}, we then have \eq{1118b} 
as asserted. $\qed$ \\

Now we estimate $\phat_x$.
By \eq{qlb},
$(1-p_y)^{n-1} \geq  e^{-1.05 n p_y}$
for $ y \in [m]$, so 
\bea
1 - (1-p_y)^{n-1} 
\leq
1 -   e^{-1.05 n p_y} \leq
 1.05 np_y,
\lbl{1114c1}
\eea
 and by \eq{qub},
 \eq{gammadef} and
the assumption  that $n \geq 1661$ by \eq{nbig1},
\bea
1 - (1-p_y)^{n-1} \geq 1 - e^{-0.9 n p_y} \geq 
(1 - e^{-0.9 \gamma}) np_y/ \gamma 
\geq
(1 - e^{-0.9}) np_y/ \gamma 
\nonumber \\
\geq  (0.55) np_y/\gamma.
      ~~~~~~~~
\lbl{1114c2}
\eea
By \eq{phatdef}, \eq{1114c1} and \eq{1114c2},
 for all $x \in [m]$ we have that
\bea
  \phat_x \leq  \frac{ 2 \gamma p_x^2 }{\sum_y p_y^2}. 
\lbl{phatub}
\eea
By \eq{gammadef}  and \eq{genvarub}, we have further that 
 \bea
\|\bphat\| :=
\sup_x
( \phat_x )\leq (2 \gamma)
 \sup_x \frac{n^2 p_x^2 }{\sum_y n^2 p_y^2} 
\leq \frac{2 \gamma^3}{n^2 \sum_y p_y^2}
\lbl{1121a}
\\
 \leq  \frac{ 
16 \gamma^3 }{\Var Y}
.
\lbl{nphatbd}
\eea
 Also, by  \eq{phatub}, 
\eq{gammadef} and \eq{genvarub},
\bea
\E \hat{p}_{X_1} \leq \frac{2 \gamma \sum_x p_x^3}{\sum_y p_y^2 } 
\leq \frac{ 2 \gamma^2}{n} ; 
\lbl{1125c}
\\
\E \hat{p}^2_{X_1} \leq 
\frac{4 \gamma^2 \sum_x p_x^5}{(\sum_y p_y^2)^2 } 
\leq
\frac{32 n^2 \gamma^2 \sum_x p_x^5}{ (\Var Y) \sum_y p_y^2 } 
\leq 
\frac{32 \gamma^5 }{n  \Var Y }. 
\lbl{1128a}
\eea

\begin{lemm}
\lbl{lemnonunidel}
With $C (\gamma)$ given by \eq{Cdef},
it is the case that
\bea
\Var (\E[Y''-Y |Y]) \leq 
\frac{(C(\gamma))^2}{\Var Y}
.
\lbl{nonunidel}
\eea
\end{lemm}
{\em Proof.}
We shall use the fact that
by \eq{genvarub} $\Var Y \leq  8 n \gamma$ so 
\bea
n^{-1} \leq 8 \gamma (\Var Y)^{-1}.
\lbl{1130b}
\eea
We estimate in turn the two variances in the right hand side of
\eq{1118b}. First consider the single sum.
Let $S_i$ denote the $i$th term in that sum, i.e.
set
\bean
S_i := 
[\tih_7(M_i,X_i)
- (1-n^{-1}) b \tih_4(M_i,X_i)  ]
 \phat_{X_i}  
\\
-
 h_0(M_i)
[ 
  2 + \phat_{X_i} \tih_{4}(M_i,X_i) 
 ]
n^{-1}, 
\eean
and set $\Sbar_i :=  S_i - \E S_i$. 
By \eq{1123a} and \eq{1125b1},
along with the fact that  $h_0(k) \in [0,2]$ so 
$0 \leq b \leq 2$,
the coefficient of $\phat_{X_i}$,
in the definition of $S_i$,
 lies
in the range
$[-1,3 ]$, 
while the coefficient of $n^{-1}$ lies in the range $[-4,0]$.
Hence,
$ | S_i + \frac{2}{n}| \leq 3  \phat_{X_i} + \frac{2}{n}.  $
By \eq{Jensen}, \eq{1128a} 
and \eq{1130b},
\bea
3n  \Var \left[ S_1 \right]     
\leq  3n \E \left[ \left(S_1 + 
(2/n) \right)^2   \right]
\leq  6n 
 \left( 
9 
 \E [\phat_{X_1}^2] + 4 n^{-2}
\right)
\nonumber \\
\leq \frac{1728 \gamma^5 + 192 \gamma}{\Var Y} .
\lbl{1118e}
\eea
Also, in the notation of Lemma  
\ref{NUmomlem},
if we write $\Sbar_i= \psi(X_i,M_i)$ we have
$\rng(\psi) \leq  4  \|\bphat\| + 4 n^{-1} $, 
and also
$\|\psi\| \leq 4  \|\bphat\| + 4 n^{-1}$.
Hence, 
  by \eq{1121a} and \eq{nphatbd},
followed  by \eq{genvarub} and  then \eq{gammadef},
\bean
\max( \rng(\psi) ,\|\psi\| )
\leq 
\left( 4 
 + \sqrt{ \frac{ \Var Y \sum_x p_x^2}{2 \gamma^6} } 
\right)
\sqrt{ \frac{32 \gamma^6}{n^2 \Var Y \sum_y p_y^2 } }
\\
\leq 
\left( 4
+ \frac{ 2n  \sum_x p_x^2}{ \gamma^3}  
\right)
\sqrt{ \frac{32 \gamma^6}{n^2 \Var Y \sum_y p_y^2 } }
\leq 
\left( 4
 + \frac{ 2}{ \gamma^2}  
\right)
\sqrt{ \frac{32 \gamma^6}{n^2 \Var Y \sum_y p_y^2 } }.
 \eean
By Lemma \ref{NUmomlem} 
 \bean
 3 n^2 \Cov \left[ S_1, S_2 \right] 
= 3 n^2 \E [\Sbar_1 \Sbar_2 ]
\leq
  9 (1 +  \gamma)   \left(
 32 \gamma^6/\Var Y  
\right)
\left(
16 + 16\gamma^{-2} + 4 \gamma^{-4}
\right)
\\
\leq \frac{4608 ( \gamma^7 + \gamma^6
 +
\gamma^5 + \gamma^4) + 1152 (\gamma^3 +  \gamma^2)
}{\Var Y } .
\eean
Combining this with  \eq{1118e} yields
\bea
   3 \Var \sum_{i=1}^n S_i & = & 3 n \Var [ S_1 ] +   3  n (n-1) \Cov [ 
S_1, S_2 ]
\nonumber \\
& \leq &
\left(
\frac{ 100 }{  
\Var Y }  
\right)
\left(
47 \gamma^7
+ 
 47 \gamma^6 + 64 \gamma^5 
+ 47 \gamma^4
+ 12 \gamma^3 + 12 \gamma^2
\right).  
\lbl{1118f}
\eea
Consider now the double sum in \eq{1118b}.
Writing $(n)_{k}$ for $n!/(n-k)!$, we have
\bea
\Var
 \sum_{(i,j): i \neq j} \xi_i  \bar{T}_j  = (n)_4  \Cov(
 \xi_1 \bar{T}_2, \xi_3\bar{T}_4)
+ (n)_2 \left(
\Var( \xi_1 \bar{T}_2 ) +
\Cov(\xi_1 \bar{T}_2, \bar{T}_1 \xi_2 )
\right)
\nonumber  \\
+  (n)_3  \left(
 \Cov(\xi_1 \bar{T}_2, \xi_1 \bar{T}_3)
 + \Cov(\xi_2 \bar{T}_1, \xi_3 \bar{T}_1)
 + 2\Cov(\xi_1 \bar{T}_2, \xi_3 \bar{T}_1) \right).
\lbl{0519a}
 \eea
For the first term of the right hand side of 
\eq{0519a}, observe that
\bean
\Cov ( \xi_1 \bar{T}_2, \xi_3 \bar{T}_4 ) = 
\E [ \xi_1 \bar{T}_2 \xi_3 \bar{T}_4 ] 
-
\E [ \xi_1 \bar{T}_2 ] \E[ \xi_3 \bar{T}_4 ] 
 \leq \E [ \xi_1 \bar{T}_2 \xi_3 \bar{T}_4 ],
\eean
and that $0 \geq \xi_i \geq - \phat_{X_i}$ by \eq{1123a}, while  
$0 \leq T_j \leq 2$.  
So by  Lemma \ref{NUmomlem}, \eq{1121a} and \eq{nphatbd},
\bea
3 n^2 \Cov ( \xi_1 \bar{T}_2, \xi_3 \bar{T}_4 )  
\leq 12 n^2   \|\bphat\|^2
(9 + 9 \gamma) \sum_x p_x^2
\leq 108 (1+ \gamma)  
\left( \frac{32 \gamma^6 }{ \Var Y} \right)
\nonumber \\
= (3456 \gamma^6 + 3456 \gamma^7)/\Var Y.
~~~~~~~
\lbl{1126a}
\eea
Now consider the last term in \eq{0519a}. 
By \eq{1128a},
\bea
3 n \Cov(\xi_1 \bar{T}_2, \xi_1 \bar{T}_3)
\leq   
3 n \E [ \xi_1^2 \bar{T}_2 \bar{T}_3 ] 
\leq
12 n  \E \phat_{X_1}^2
\leq
\frac{ 384  \gamma^5 }{  \Var Y},
\label{1126b}
\eea
while by \eq{1125c} and \eq{1130b},
\bea
 3  n (\Cov(\xi_2 \bar{T}_1, \xi_3 \bar{T}_1)
 + 2\Cov(\xi_1 \bar{T}_2, \xi_3 \bar{T}_1) )
\leq
  3 n \E \bar{T}_1^2 \xi_2 \xi_3  +
 6 n   \E
 \xi_1 \bar{T}_2 \xi_3 \bar{T}_1 
 \nonumber \\
\leq
36 n 
\E [   \phat_{X_1} \phat_{X_2} ]
\leq 
144 \gamma^4
n^{-1}
\leq \frac{1152 \gamma^5}{\Var Y} .
~~~~~~~
\lbl{1128b}
\eea
The middle term in the right side of \eq{0519a} 
is smaller; since $\E [\phat_{X_1} \phat_{X_2} ] \leq \E [\phat_{X_1}^2 ]$,
 \eq{1128a} gives
\bean
3  ( \Var( \xi_1 \bar{T}_2 ) + \Cov(\xi_1 \bar{T}_2, \bar{T}_1 \xi_2 ) )
\leq 24 \E [ \phat_{X_1}^2]
\leq 
\frac{ 768 \gamma^5 }{n \Var Y } ,
\eean
and  since $n \geq 1661 $ by \eq{nbig1}, 
combined with \eq{0519a}, 
\eq{1126a}, \eq{1126b}, and \eq{1128b}, this shows that 
\bea
3 n^{-2}
\Var
 \sum_{(i,j): i \neq j}
 \xi_i  \bar{T}_j  
\leq
\left( \frac{100  }{ \Var Y} \right)
( 35 \gamma^7 + 35 \gamma^6 + (15.9) \gamma^5). 
\lbl{1126c}
\eea
Also, by \eq{1130b} we obtain 
$$
12 (n-1)^2  \leq \left( \frac{12}{1660}
\right) \left( \frac{n}{n-1} \right) n^{-1} \leq
\left( \frac{96\times 1661}{1660^2 } \right)  
\frac{\gamma}{\Var Y} \leq \frac{\gamma}{10 \Var Y} .
$$
Combining this with \eq{1118b}, \eq{1118f} and \eq{1126c} yields
\bean
\Var (\E[Y''-Y |Y]) \leq 
 \left( \frac{100 }{\Var Y} \right) 
(82 (\gamma^7 +  \gamma^6) + 80 \gamma^5 + 47 \gamma^4 + 12 ( \gamma^3 +   \gamma^2 )) .  
\quad
\qed
\eean
 \vspace{0.5em}

{\em Proof of Theorem \ref{nonunithm}.}
It remains to prove \eq{thm2eq}.
By \eq{1114c1},
\bean
\E Y = n \sum_x p_x (1 - (1-p_x)^{n-1} ) 
\leq 1.05 n^2 \sum_x p_x^2,  
\eean
so that by \eq{genvarlb},
\bea
\mu_Y / \sigma_Y^2
\leq
  8165 \gamma^2 e^{2.1 \gamma}.
\lbl{EVrat}
\eea
We can apply Lemma \ref{LGthm} with $B=3$ by \eq{1114b}.
According to that lemma,  for \eq{delbd} 
we need $B \leq \sigma^{3/2} / \sqrt{6 \mu}$, i.e.
$\sigma^3 \geq 54 \mu$.
By \eq{EVrat}, a sufficient condition for this is that
$\sigma \geq 54(8165) \gamma^2 e^{2.1 \gamma}$. 
But if this
condition fails then since $577 > 54$, the right hand side of 
\eq{thm2eq} is greater than 1 so we are anyway guaranteed that
\eq{thm2eq} holds.

 Therefore, from now on we may assume that
$\sigma_Y \geq 54(8165) \gamma^2 e^{2.1 \gamma}$, so that
\eq{delbd} is valid.
Using \eq{EVrat} and Lemma \ref{lemnonunidel}, this gives
\bean
D_Y
\leq  \sigma_Y^{-1}  
\left( 1.2 + 8165 \gamma^2 e^{2.1 \gamma} ( 576 + 108 \sigma_Y^{-1} + 23C) 
\right)
\eean
with $C = C(\gamma)$.
We may then deduce \eq{thm2eq}.
$\qed$ \\


 {\bf Acknowledgement.} It is a pleasure to thank the Institut f\"ur
Stochastik, Universit\"at Karlsruhe,  for their hospitality during
the preparation of this manuscript, and  Norbert Henze for valuable
comments on an earlier draft.



\end{document}